\documentclass[]{amsart}

\usepackage{lineno,hyperref}
\usepackage{amsmath}
\usepackage{amsfonts}
\usepackage{latexsym}
\usepackage{amsthm}
\usepackage{listings}
\usepackage{amssymb, ifthen, enumerate}
\usepackage{color}
\usepackage{ulem}
\usepackage{longtable}
\modulolinenumbers[5]

\newtheorem{lemma}{Lemma}[section]
\newtheorem{thm}[lemma]{Theorem}

\newtheorem{definition}[lemma]{Definition}

\renewcommand{\qed}{\hspace*{\fill}\rule{1.5ex}{1.5ex}}
\newcommand{\N}{\mathbb{N}}	
\newcommand{\Ri}{\mathcal{R}}     
\newcommand{\Alt}{\textrm{Alt}}     

\newcommand{\Menge}[2]{\left\{ #1 \;\middle|\; #2 \right\}}

\DeclareMathOperator{\Aut}{Aut}

\DeclareMathOperator{\SL}{SL}

\DeclareMathOperator{\SU}{SU}
\DeclareMathOperator{\M}{M}
\DeclareMathOperator{\PSL}{PSL}
\DeclareMathOperator{\PSU}{PSU}

\DeclareMathOperator{\PGL}{PGL}
\DeclareMathOperator{\PGU}{PGU}

\makeatletter 

\newcommand{\abs}[1]{\left| #1 \right|}

\newcommand{\prf}[1]{
\noindent Proof.
#1
\qed\\[0.5eM]
}

\begin{document}

\begin{center}
\Large{\textbf{Finite simple groups acting with fixity $3$ and their occurrence as groups of automorphisms of Riemann surfaces (extended version)}}

\vspace{0.2cm} \small{Patrick Salfeld and Rebecca Waldecker}
\end{center}





\normalsize
\vspace{0.5cm}
\textbf{Abstract.}

Motivated by the theory of Riemann surfaces, we classify all possibilities for finite simple groups acting faithfully on a compact Riemann surface of genus at least $2$ in such a way that all non-trivial elements have at most three fixed points on each non-regular orbit and at most four fixed points in total.
In each case we also give information about the branching datum of the surface.
There is a shorter version of this article (submitted for publication), and in this extended version we give many more details about the GAP code that we use for the calculations.
We also explicitly include a lemma that we only quote in the short version, so we can explain how exactly the GAP calculations and the lemma work together.

\section{Introduction}

In previous work we have studied finite groups from the perspective
of low fixity actions. In this context we say that a group $G$ \textbf{acts with fixity $k \in \N$
on a set} $\Omega$ if and only if $k$ is the maximum number of fixed points of elements of $G^\#$ 
(i.e. the set of non-identity elements of $G$).
While these questions are interesting in their own right from the perspective of permutation group theory,
they have originally been motivated from
the theory of Riemann surfaces.
Magaard and V\"olklein explain in \cite{MV} how a group-theoretical property about fixed points of non-trivial automorphisms of a curve (i.e. a compact Riemann surface of genus at least 2) can be used for the construction of Weierstrass points of this curve. Here, by an automorphism of a Riemann surface we mean a biholomorphic bijective map of the surface to itself.

The connection between fixed points of automorphisms and Weierstrass points goes back to Schoeneberg, and we state his result along the lines of \cite{MV}: If a non-trivial automorphism of a curve fixes at least five points on the curve, then these fixed points are Weierstrass points. 
While \cite{MV}, \cite{Wa} and \cite{Ro} give examples for arguments based on Schoeneberg's result, their work raises the question of what happens if there are no automorphisms with enough fixed points, making this method unapplicable.  
This is where Kay Magaard and the second author started their work, and it is the reason why
 the classification results in \cite{MW}, for fixity 2, were followed by investigating
which ones of the finite simple groups from \cite{MW} actually occur as groups of automorphisms of a curve, 
and with low fixity in their action (\cite{SW}).

In parallel to \cite{SW}, we now build on 
the classification results in \cite{MW3}, for fixity 3, with the aim to
give information about the curves such that
the finite simple groups from \cite{MW3} act as groups of automorphisms and in such a way that every non-trivial element fixes at most three points on each orbit and at most four points in total. We describe the curves and the action of the group
in terms of branching data, which does not give an explicit isomorphism type of the curve.

In subsequent work we will complete this part of the project by including all finite simple groups that act with fixity 4,
completing our overview over situations where Schoeneberg's argument is not applicable. 
After this background and motivation, we intend to keep this article short. We do not repeat all the background that we gave in \cite{MW} and \cite{SW}, in particular we do not go into details about the relevant theory of Riemann surfaces here. 
Instead, we go straight to our main result:

\begin{thm}\label{main3fpH}
Suppose that $G$ is a finite simple group that acts faithfully on a compact Riemann surface
$X$ of genus at least $2$.
Suppose further that there is at least one orbit on which $G$ acts with fixity 3, that $G$ acts with fixity at most 3 on each orbit of $X$ and with fixity at most 4 in total.
Then $G$ acts with one of the branching data given in Table \ref{DataList}.

Conversely, if $l$ is a list from Table \ref{DataList} and $l$ is a Hurwitz datum, then there exists a
compact Riemann surface $X$
of genus at least $2$ on which the group $G$ in $l$ acts faithfully, and for all choices of $X$, it is true that
$G$ acts with fixity at most 3 on each orbit of $X$, with fixity 3 on at least one orbit and with fixity at most 4 in total.
\end{thm}

We now explain how to interpret the data in Table \ref{DataList}, which belongs to the main theorem, by providing an illustrative example:

The list $[\PSL_2(7),g,g_0 \mid [4,2],[7,1]]$ is a potential branching datum for
the group $\PSL_2(7)$, acting as a group of automorphisms on a compact Riemann surface 
with genus $g$ and cogenus $g_0$ and such that there are exactly three non-regular orbits. $\PSL_2(7)$
acts on two of those with point stabilisers of order 4 and on one of them with point stabilisers of order 7.
For such a list to be a Hurwitz datum (as mentioned in the theorem), the Hurwitz formula from Definition 2.1 must be satisfied.

\begin{longtable}{llp{5.8cm}}
	
	Line&Hurwitz datum&Remark\\[0.4ex]\hline\hline
	\\[-2ex]
		1&$[\Alt_7,g,g_0 \mid [ 7,1]]$&$g_0 \ge 1$\\
		2&$[\Alt_8,g,g_0\mid [7,1]]$&$g_0  \ge 1$\\
		3&$[\M_{22},g,g_0\mid[7,1]]$&$g_0  \ge 1$\\
		4&$[\PSL_4(3),g,g_0 \mid [{13},1]]$&$g_0  \ge 1$\\
		5&$[\PSL_4(5),g,g_0 \mid [{31},1]]$&$g_0  \ge 1$\\
		6&$[\PSU_4(3),g,g_0 \mid [7,1]]$&$g_0  \ge 1$\\
		7&$[\PSL_3^\epsilon(q),g,g_0 \mid [\alpha,1]]$&$g_0 \ge 1$, $\epsilon \in\{1,-1\}$, $q \ge 3$ prime power, $\alpha=\frac{q^2+\epsilon q+1}{(3,q-\epsilon)}$.\\
		8&$[\PSL_2(7),g,g_0 \mid [7,1]]$&$g_0  \ge 1$\\

		9&$[\PSL_2(7),g,g_0 \mid [3,1],[7,1]]$&$g_0  \ge 1$\\
10&$[\PSL_2(7),g,g_0 \mid [3,2],[7,1]]$&$g_0  \ge 0$\\
                11&$[\PSL_2(7),g,g_0 \mid [4,1],[7,1]]$&$g_0  \ge 1$\\
12&$[\PSL_2(7),g,g_0 \mid [4,2],[7,1]]$&$g_0  \ge 0$\\

13&$[\PSL_2(7),g,g_0 \mid [3,1],[4,1],[7,1]]$&$g_0  \ge 0$\\
14&$[\PSL_2(7),g,g_0 \mid [3,2],[4,1],[7,1]]$&$g_0  \ge 0$\\
15&$[\PSL_2(7),g,g_0 \mid [3,1],[4,2],[7,1]]$&$g_0  \ge 0$\\
16&$[\PSL_2(7),g,g_0 \mid [3,2],[4,2],[7,1]]$&$g_0  \ge 0$\\	
17&$[\PSL_3(4),g,g_0 \mid [5,1],[7,1]]$&$g_0  \ge 1$\\
18&$[\PSL_3(4),g,g_0 \mid [5,2],[7,1]]$&$g_0  \ge 0$\\	
		\\[-2ex]\hline\hline\\[-2ex]

	\caption[List of Hurwitz data derived from Theorem ...]{\small Hurwitz data for group actions with only one non-regular orbit and fixity 3 or with up to five non-regular orbits and mixed fixity, $g \ge 2$.}
	\label{DataList}
	
\end{longtable}

For the proof of the theorem,
we recall some relevant definitions, then we refer to \cite{MW3} for the specific groups to consider and analyse them 
with the methods developed in \cite{SW} and by using \cite{GAP4}.

\textbf{Acknowledgements.} This article was written in memory of Kay Magaard. We remember him with much gratitude for
initiating the project on groups acting with low fixity, for
suggesting questions along the lines of this article and for many discussions on
the subject. We also thank Chris Parker for very helpful remarks on earlier drafts of this paper.


\section{Preliminaries}

Most of our notation is standard, we refer to \cite{Miranda1997} for background information on Riemann surfaces
and we only mention the following for clarity: \\
If $G$ is a group and $g,h \in G$, then we write \textbf{conjugation} from the right, so $g^h:=h^{-1}  g h$, and we write \textbf{commutators} as $[g,h]:=g^{-1}  h^{-1}  g  h$.\\
Throughout, we suppose that $X$ is a compact Riemann surface of genus $g\ge 2$ and that $G\le\Aut(X)$. Let $X/G=\Menge{ x^G}{ x\in  X}$ denote the \textbf{space of $G$-orbits on $X$} and let $ g_0\ge 1$ denote the \textbf{cogenus}, i.e. the genus of $X/G$.
Then all elements in $X/G$ are finite and their cardinalities (which we refer to as \textbf{orbit lengths}) divide $|G|$.
We recall that an orbit is \textbf{non-regular }if and only if some element of $G^\#$ fixes a point in it, which happens if and only if its length is strictly less than $|G|$. For all $x\in X$ we denote the point stabiliser of $x$ in $G$ by $G_x$.
In \cite{SW} we give more details about the notation. Here we just recall:
	\begin{itemize}
                \item $G$ has infinitely many regular orbits on $X$ and only finitely many non-regular orbits. (Proposition III.3.2, Theorem III 3.4 and 
Corollary II.1.34 in \cite{Miranda1997}.)
		\item  $\abs {\Aut( X)}$ is finite and for all $x \in X$, the subgroup $G_x$ is cyclic. (Theorem VII.4.18 and Proposition
III.3.1 in \cite{Miranda1997}.)
	\end{itemize}
In particular, all groups considered here are finite.

\begin{definition}\label{einf1-not:stlit}

Given a finite group $G$, non-negative integers $g,g_0,r,m_1,n_1,...,$\\
$m_r,n_r$ and a list
$l:=[G, g, g_0 \mid [m_1,n_1],\ldots, [m_r,n_r]]$,
we refer to $l$ as a \textbf{Hurwitz datum} if and only if the Hurwitz formula 
is satisfied:
$$2(g-1)=\abs G\left( 2(g_0-1)+\sum_{j=1}^r n_j\left(1-\frac 1{m_j}\right)\right).$$

We say that $G$ \textbf{acts with branching datum} $l$ on a Riemann surface $X$ if and only if
$X$ has genus $g$, $X/G$ has genus $g_0$ and for each $i \in \{1,...,r\}$ there are exactly $n_i$ non-regular orbits
of $X/G$
on which $G$ acts with point stabilisers of order $m_i$.\\
In such a case we refer to $X$ as a \textbf{witness for $l$}.\\
We let $\Ri_4$ denote the set of Hurwitz data $l=[G, g, g_0 \mid [m_1,n_1],\ldots, [m_r,n_r]]$ where
$m_1<m_2< \cdots <m_r$ and such that $l$ has a witness $X$ on which $G$
 acts with fixity at most $4$. Finally, we denote by $\Ri_4^*$ the subset of $\Ri_4$ of lists such that, for all witnesses $X$ for $G$ with respect to $l$, $G$ acts with fixity at most $4$ on $X$.
Whenever we introduce a Hurwitz datum or a list from $\Ri_4$, then we use all the notation explained here.
\end{definition}

The main lemma that we use later when checking specific groups with GAP is the following
(originally from \cite{Broughton1990}, see also Lemma 3.2 in
\cite{SW}):

\begin{lemma}\label{crit}
Suppose that $l=[G, g, g_0 \mid [m_1,n_1],\ldots, [m_r,n_r]]$ is a Hurwitz datum.
Then $G$ has a witness $X$ with respect to $l$ if and only if there exist
elements
 $a_1,\ldots,a_{g_0},b_1,\ldots,b_{g_0},c_{1,1},\ldots,c_{1,n_1},\ldots,c_{r,1},
\ldots,c_{r,n_r}\in G$ satisfying the following conditions:

\begin{enumerate}
\item
For all $j\in\{1,\ldots,r\}$ and all $i\in\{1,\ldots,n_j\}$ it is true that $o(c_{j,i})=m_j$,
    \item
    $\prod_{k=1}^{g_0}[a_k,b_k]\cdot \prod_{i=1}^{n_1}c_{1,i}\cdots \prod_{i=1}^{n_r}c_{r,i}=1$, and
    \item
    $\langle a_1,\ldots, a_{g_0},b_1,\ldots,b_{g_0},c_{1,1},\ldots,c_{r,n_r}\rangle=G$. \end{enumerate}

Furthermore, if $G$ has a witness $X$ with respect to $l$ and if $h\in G^\#$ fixes a point in $X$, then $h$ is conjugate to a power of one of $c_{1,1},\ldots,c_{r,n_r}$ in $G$.
\end{lemma}

We point to a special case of this lemma, namely when $g_0=0$. Then the elements
$a_1,\ldots,a_{g_0},b_1,\ldots,b_{g_0}$ do not exist and the properties (a) -- (c) just refer to a generating set
$\{c_{1,n_1},\ldots,c_{r,1},
\ldots,c_{r,n_r}\}$ of $G$. In particular (b) simplifies substantially in this special case.\\

We close this section with the list of groups that we will analyse and with a discussion of the cases with ``mixed fixity'':

\begin{thm}\label{list3}
Suppose that $G$ is a non-abelian finite simple group that acts faithfully, transitively and with fixity $3$ on a set $\Omega$, with cyclic point stabilisers.
If $\alpha \in \Omega$, then
$G$, $\abs\Omega$ and $\abs{G_\alpha}$ are as in the following list:

(1)		$\Alt_7$, $|\Omega|=2^3\cdot 3^2\cdot 5$			and	$\abs{G_\alpha}=7$.	

(2)		$\Alt_8$,		$|\Omega|=2^6\cdot 3^2\cdot 5$	and			$\abs{G_\alpha}=7$.					
	
(3)		$\M_{22}$,		$|\Omega|=2^7\cdot 3^2\cdot 5\cdot 11$	and	$\abs{G_\alpha}=7$.					
	
(4)		$\PSL_4(3)$,		$|\Omega|=2^7\cdot 3^6\cdot 5$		and		$\abs{G_\alpha}=13$.					
	
(5)		$\PSL_4(5)$,		$|\Omega|=2^7\cdot 3^2\cdot 5^6\cdot 13$	and	$\abs{G_\alpha}=31$.					
	
(6)		$\PSU_4(3)$,		$|\Omega|=2^7\cdot 3^6\cdot 5$			and	$\abs{G_\alpha}=7$.					

(7)		$\PSL_3(q)$,		$|\Omega|=q^3(q-1)(q^2-1)$,	$\abs{G_\alpha}=(q^2+q+1)/d$, where $q$ is a prime power, $q\ge 2$, and $d:=(3,q-1)$.
	
(8)		$\PSU_3(q)$,		$|\Omega|=q^3(q+1)(q^2-1)$,			$\abs{G_\alpha}=(q^2-q+1)/d$, where $q$ is a prime power, $q\ge 3$, and $d:=(3,q+1)$.

\end{thm}

\prf{By hypothesis $G$ acts with fixity $3$, so we may apply Theorem 1.1 in \cite{MW3} with the additional information that the point stabilisers are cyclic to obtain the result.
}

There are two groups in our list (based on the previous theorem) that also allow for a fixity 2 action, so there will be several non-regular orbits and the group acts with fixity 2 or 3 (which we will refer to as \textbf{mixed fixity action}).

\begin{lemma}\label{mixed}
Suppose that $G$ is a finite simple group that exhibits transitive faithful actions with fixity 2 and 3, and with cyclic point stabilisers, respectively. Then one of the following holds:

$G \cong \PSL_2(7)$ acts with fixity 2 and cyclic point stabilisers of order 3 or 4, and with fixity 3 and cyclic point stabilisers of 
order 7, or

$G \cong \PSL_3(4)$ acts with fixity 2 and point stabilisers of order 5, and with fixity 3 and point stabilisers of order 7.
\end{lemma}

\prf{
We combine the results from \cite{MW} and \cite{MW3} with the additional hypothesis that the point stabilisers are cyclic. More specifically we inspect Theorem 1.1 in \cite{MW3} and Lemmas 3.11 and 3.13 of \cite{MW}. This gives exactly the groups that are mentioned in the lemma. 
}

A specific mixed fixity example is $G:=\PSL_2(7)$, acting with branching datum $[G,24,0 \mid [3,1], [4,1], [7,1]]$ in our notation.
In detail this means that the compact Riemann surface $X$ has genus 24, that $X/G$ has genus 0 and that $G$ has three non-regular orbits in its action on $X$:
The first one has size 56 and the group acts with fixity 2, the second one has size 42 and the group acts with fixity 2, and the last one has size 24 and the group acts with fixity 3.\\
In particular, all non-trivial group elements have less than five fixed points on $X$.

In general, the hypothesis of global fixity at most 4 does not allow for too many non-regular orbits, which is why in the following section, whenever $[G, g, g_0 \mid [m_1,n_1],\ldots, [m_r,n_r]]$ is a Hurwitz datum, then $r \le 3$ and $n_1,n_2,n_3 \le 2$.


\section{The \texttt{GAP} calculations}



In order to prove Theorem \ref{main3fpH}, we look at all cases from Table \ref{DataList}, first in the case of minimal
cogenus $g_0$. This can mostly be done by calculations in GAP (\cite{GAP4}), which is the content of this section, and we use the \texttt{MapClass} package (\cite{JMSV18}) in some cases.
More precisely, we have written three \texttt{GAP} functions that we will explain below, and at the end of the section we will display the output of these functions, applied to the corresponding groups in the same order as they appear in Table \ref{DataList}.
All functions use
the character table of the group, and \texttt{GenTuple1} also uses the functionality of the \texttt{MapClass} Package (see \cite{JMSV18}), specifically the function
\texttt{GeneratingMCOrbits}.
The objective is always the same: Finding a set of generators of the group that satisfies the
conditions in Lemma \ref{crit} for the minimal cogenus $g_0$ as given in Table \ref{DataList}.\\


\noindent
\underline{The function \texttt{GenTuple1}}

\smallskip
We suppose that the package \texttt{MapClass} is already loaded.

\begin{verbatim}
GenTuple1:=function(G,cogen,orders)
   local ct,cc,list,pos,ro,tuple,orb,orb1,h,i,len;
   ct := CharacterTable(G);;
   cc := ConjugacyClasses(ct);;
   list := OrdersClassRepresentatives(ct);;
   pos := List(orders, i-> Positions(list,i) );;
   ro := List(pos, i-> Random(i) );;
   tuple := List(ro,i->Representative(cc[i]) );;
   
   orb := GeneratingMCOrbits(G,cogen,tuple : 
   OutputStyle:="single-line");;
   len := Length(orb);;
   if len=0 then return; fi;

   for i in [1..len] do
      Print(Length(orb[i].TupleTable));
      Print(" tuples in orbit ");
      Print(i);
      Print("\n");
   od;

   Print("\n");
   Print("Picking random tuple in random orbit ...\n");
   orb1 := orb[Random([1..len])];;
   len := Length(orb1.TupleTable);;
   h := SelectTuple(orb1, Random([1..len]) );;

   Print("Testing tuple ... ");
   if Comm(h[1],h[2])*Product(h{[3..Length(h)]})=One(G)
    or Product(h)=One(G) then
      Print("Lemma 2.2 (b) and (c) are satisfied!\n");
	  return h;;
   else
      Print("Something went wrong. Please try again!");
	  return;;
   fi;
end;
\end{verbatim}

Some comments:

As input we have a group, a cogenus and a list of element orders. 
For each element order in the list "orders", we find out where they are located in the list of orders of the conjugacy class representatives,
and then we randomly pick elements, one of each given order. This gives "tuple", which is one of the ingredients for the function
 \texttt{GeneratingMCOrbits}.
If this function does not find any orbit, then we stop. Otherwise we print the number of generating tuples per orbit
and we randomly pick one. Since we alrady took care of the element orders in the input (Lemma \ref{crit}(a)) and
all tuples that we find generate the group (according to the documentation of the \texttt{MapClass} package), we only
check statement (b) from Lemma \ref{crit}. This is a commutator or a product, depending on whether 
the cogenus is 0 or 1. 
The groups $\Alt_7$ and $\PSL_2(7)$ can all be handled with this function. However, some problems arise for larger groups, starting with $\Alt_8$.

\medskip
\noindent
\underline{The function \texttt{GenTuple2}}

\smallskip

The first list which is a bit too much for the previous function is $[\Alt_8, g, 1 \mid [7, 1]]$.
The bottleneck is 
\texttt{GeneratingMCOrbits} with genus 1, which with a tuple of length 1 (as we did for $\Alt_7$) takes several minutes.
Alternatively, we can use \texttt{GeneratingMCOrbits} with genus 0
and in this way find a suitable set of generators according to Lemma \ref{crit}.

We suggest another alternative here, which is an explicit construction of a suitable generating tuple. 
It works in basically the same way for all the groups in Lines 2 -- 6 of Table \ref{DataList}, which have in common
that we always have minimal cogenus 1
and only one non-regular orbit.
The little technical details are dealt with in a case distinction.

\begin{verbatim}
GenTuple2:=function(G,order)
   local ct,cc,i,id,k,cck,a,b,c,h,pos;
   id:=Order(G);;
   if id=443520 and order=7 then #M22
      k:=11;;
   elif id=3265920 and order=7 then #U4(3)
      k:=9;;
   elif id=6065280 and order=13 then #L4(3)
      k:=5;;
   elif id=7254000000 and order=31 then #L4(5)
      k:=13;;
   elif id=20160 and order=7 then #A8
      k:=15;;
   else
      ErrorNoReturn("You've picked the wrong group");
   fi;

   ct := CharacterTable(G);;
   cc := ConjugacyClasses(ct);;

   Print("Computing tuple ...\n");
   pos:=Random(Positions(OrdersClassRepresentatives(ct),order));;
   c:=Representative(cc[pos]);;
   pos:=Random(Positions(OrdersClassRepresentatives(ct),k));;
   for i in cc[pos] do
      if i*c^(-1) in cc[pos] and IsConjugate(G,i,i*c^(-1)) then
         a:=i;;
         h:=i*c^(-1);;
         break;
      fi;
   od;
   b:=RepresentativeAction(G,a,h);;

   Print("Testing tuple ... ");
   if Comm(a,b)*c=One(G) and Group([a,b,c])=G then
      Print("Lemma 2.2 (b) and (c) are satisfied!\n");
	  return [a,b,c];;
   else 
      Print("Something went wrong. Please try again!");
	  return;;
   fi; 
end;
\end{verbatim}

Again, some comments:

This time the input consists of a group and an element order.
Another element order that will be needed in the process 
is chosen in the beginning, with a case distinction.
After this prepartion, we proceed as before and randomly choose elements of 
suitable orders, performing the necessary checks for Lemma \ref{crit}.
The output is a generating triple for the group that satisfies all conditions in the lemma,
provided that the input order we started with is the order of the point stabiliser from the branching datum.

\noindent
\newpage
\underline{The function for $\PSL_3(4)$:}

\smallskip
It turned out that it was easier to deal with the group $\PSL_3(4)$ and its two appearances in Table \ref{DataList}
separately rather than trying to extend previous functions.

\begin{verbatim}
GenTupleL34:=function(cogen)
   local G,i,cc5,cc3,a,b,c,h,x,y;
   G:=PSL(3,4);;
   Print("Computing tuple ...\n");
   c:=Random(Filtered(G,i->Order(i)=7));;
   cc5:=Filtered(G,i->Order(i)=5);;
   for i in cc5 do
      if i*c in cc5 then
         b:=i;;
         a:=(i*c)^(-1);;
         break;
      fi;
   od;

   if cogen=1 then
      cc3:=Filtered(G,i->Order(i)=3);;
      for i in cc3 do
         if i*a in cc3 and IsConjugate(G,i,i*a) then
            x:=i;;
            h:=i*a;;
            break;
         fi;
      od;
      y:=RepresentativeAction(G,x,h);;
   fi;

   Print("Testing tuple ... ");
   if a*b*c=One(G) and Group([a,b,c])=G then
      Print("Lemma 2.2 (b) and (c) are satisfied!\n");
   else
      Print("Something went wrong. Please try again!");
	  return;;
   fi;

   if cogen=0 then
      return [a,b,c];;
   elif cogen=1 then
      return [x,y,b,c];;
   fi;
end;
\end{verbatim}

The basic setup is as before, with more details specific to the group $\PSL_3(4)$.
In the input we only distinguish between the cases with
minimal cogenus 0 or 1, and  this is taken into account for the output as well.
We obtain a generating triple for cogenus 0 and a generating quadruple for cogenus 1, again with all conditions from Lemma \ref{crit} satisfied.

\bigskip

As promised, we close with the output using these functions, for almost all groups. We write down the input,
but omit the output
for Lines 5 and 6 of the table because the group elements are very long in their cycle presentation.
The first input for the list $[\PSL_2(7),g,0 \mid [3,2],[7,1]]$, together with its output, shows that we need to enter the element orders
in the tuple "orders" of the function \texttt{GenTuple1} with multiplicities!

\begin{verbatim}
gap> GenTuple1(AlternatingGroup(7),1,[7]);
 Computation complete : 5 orbits found.
84 tuples in orbit 1
84 tuples in orbit 2
21 tuples in orbit 3
56 tuples in orbit 4
21 tuples in orbit 5

Picking random tuple in random orbit ...
Testing tuple ... Lemma 2.2 (b) and (c) are satisfied!
[ (1,7,2,6)(3,5), (3,4,6,7,5), (1,2,3,4,5,7,6) ]

gap> GenTuple1(PSL(2,7),1,[7]);
 Computation complete : 1 orbits found.
7 tuples in orbit 1

Picking random tuple in random orbit ...
Testing tuple ... Lemma 2.2 (b) and (c) are satisfied!
[ (1,2,5)(3,7,8), (1,2,6,8)(3,7,5,4), (2,4,6,5,8,3,7) ]

gap> GenTuple1(PSL(2,7),1,[3,7]);
 Computation complete : 2 orbits found.
693 tuples in orbit 1
630 tuples in orbit 2

Picking random tuple in random orbit ...
Testing tuple ... Lemma 2.2 (b) and (c) are satisfied!
[ (1,8,5)(2,6,3), (1,2)(3,6)(4,5)(7,8), 
(1,2,5)(3,7,8), (2,3,5,4,7,8,6) ]

gap> GenTuple1(PSL(2,7),0,[3,7]);
 Computation complete : 0 orbits found.

gap> GenTuple1(PSL(2,7),0,[3,3,7]);
 Computation complete : 1 orbits found.
1 tuples in orbit 1

Picking random tuple in random orbit ...
Testing tuple ... Lemma 2.2 (b) and (c) are satisfied!
[ (1,8,4)(2,7,3), (1,2,8)(4,5,6), (2,4,6,5,8,3,7) ]

gap> GenTuple1(PSL(2,7),1,[4,7]);
 Computation complete : 2 orbits found.
448 tuples in orbit 1
448 tuples in orbit 2

Picking random tuple in random orbit ...
Testing tuple ... Lemma 2.2 (b) and (c) are satisfied!
[ (1,6)(2,5)(3,7)(4,8), (1,8,6)(2,7,5), 
(1,2,8,4)(3,7,6,5), (2,4,6,5,8,3,7) ]

gap> GenTuple1(PSL(2,7),0,[4,4,7]);
 Computation complete : 1 orbits found.
1 tuples in orbit 1

Picking random tuple in random orbit ...
Testing tuple ... Lemma 2.2 (b) and (c) are satisfied!
[ (1,7,6,3)(2,8,5,4), (1,2,5,7)(3,8,6,4), (2,3,5,4,7,8,6) ]

gap> GenTuple1(PSL(2,7),0,[3,4,7]);
 Computation complete : 2 orbits found.
1 tuples in orbit 1
1 tuples in orbit 2

Picking random tuple in random orbit ...
Testing tuple ... Lemma 2.2 (b) and (c) are satisfied!
[ (1,8,3)(4,7,5), (1,2,6,8)(3,7,5,4), (2,3,5,4,7,8,6) ]

gap> GenTuple1(PSL(2,7),0,[3,3,4,7]);
 Computation complete : 2 orbits found.
56 tuples in orbit 1
56 tuples in orbit 2

Picking random tuple in random orbit ...
Testing tuple ... Lemma 2.2 (b) and (c) are satisfied!
[ (1,3,8)(4,5,7), (1,2,5)(3,7,8), 
(1,2,5,7)(3,8,6,4), (2,4,6,5,8,3,7) ]

gap> GenTuple1(PSL(2,7),0,[3,4,4,7]);
 Computation complete : 2 orbits found.
42 tuples in orbit 1
42 tuples in orbit 2

Picking random tuple in random orbit ...
Testing tuple ... Lemma 2.2 (b) and (c) are satisfied!
[ (2,4,5)(3,8,6), (1,7,4,5)(2,8,3,6), 
(1,2,5,7)(3,8,6,4), (2,4,6,5,8,3,7) ]

gap> GenTuple1(PSL(2,7),0,[3,3,4,4,7]);
 Computation complete : 2 orbits found.
2352 tuples in orbit 1
2352 tuples in orbit 2

Picking random tuple in random orbit ...
Testing tuple ... Lemma 2.2 (b) and (c) are satisfied!
[ (2,4,3)(5,7,8), (1,5,2)(3,8,7), 
(1,2,6,8)(3,7,5,4), (1,2,8,4)(3,7,6,5),
  (2,4,6,5,8,3,7) ]

gap> GenTuple2(AlternatingGroup(8),7);
Testing tuple ... Lemma 2.2 (b) and (c) are satisfied!
[ (1,3,2,5,7)(4,8,6), (2,4,5,6,8,7,3), (1,2,3,4,5,6,8) ]

gap> GenTuple2(MathieuGroup(22),7);
Testing tuple ... Lemma 2.2 (b) and (c) are satisfied!
[ (1,3,18,14,6,10,17,21,13,5,11)(2,7,8,16,12,4,9,20,22,15,19),
  (1,19,12,13,3,11)(2,5,16,9,7,14)(4,8,10)(6,17)(15,22,20)(18,21),
  (1,12,15,8,5,4,21)(2,13,20,19,3,14,11)(6,9,17,7,16,18,10) ]

gap> GenTuple2(PSL(4,3),13);
Testing tuple ... Lemma 2.2 (b) and (c) are satisfied!
[ (1,19,5,30,3)(2,18,37,39,25)(4,17,28,21,26)(6,10,9,29,20)
(7,32,13,31,15)(8,33,38,40,11)(12,34,24,22,36)(14,35,23,16,27), 
(1,36,27,10,15,32,17,34,28,20,39,19,4)(2,35,5,16,40,6,29,9,11,13,31,23,21)
(3,37,18,38,26,25,8,30,22,33,7,12,14),
(1,2,6,4,12,13,10,11,7,3,9,8,5)(15,20,24,21,30,31,28,35,33,22,27,38,23)
(16,17,34,19,37,36,39,29,25,18,40,26,32) ]

gap> GenTuple2(PSL(4,5),31);
Testing tuple ... Lemma 2.2 (b) and (c) are satisfied!
...

gap> GenTuple2(PSU(4,3),7);
Testing tuple ... Lemma 2.2 (b) and (c) are satisfied!
...

gap> GenTupleL34(1);
Computing tuple ...
Testing tuple ... Lemma 2.2 (b) and (c) are satisfied!
[ (3,5,4)(6,20,11)(7,18,12)(8,19,10)(9,21,13)(14,16,15),
  (1,13,3,5,20,10,19)(2,15,17,9,18,11,14)(4,6,7,16,21,12,8),
  (2,13,20,7,14)(3,12,18,9,15)(4,11,19,8,16)(5,10,21,6,17),
  (1,3,14,16,6,2,8)(4,19,15,21,18,12,7)(5,13,17,11,10,20,9) ]

gap> GenTupleL34(0);
Computing tuple ...
Testing tuple ... Lemma 2.2 (b) and (c) are satisfied!
[ (1,21,18,19,20)(2,9,12,6,5)(3,17,15,4,13)(7,14,10,16,8),
  (2,13,15,8,18)(3,12,16,7,19)(4,11,14,9,20)(5,10,17,6,21),
  (1,20,2,21,12,13,5)(3,18,16,9,7,10,6)(4,19,8,17,14,11,15) ]

\end{verbatim}

Done! All the lists with individual groups from Table \ref{DataList} have been considered, at least for the minimal possible cogenus.
All that is left to do, in the next section, is to deal with Line 7 of the table and then to collect all the information together for the proof of the main result.


\section{The proof of the main thorem}

After the \texttt{GAP} calculations, it only remains to consider the groups $\PSL_3(q)$ and $\PSU_3(q)$, where
we need a generic argument.

\begin{lemma}\label{RHD:LU3}~\\
	Let $\epsilon \in\{1,-1\}$, let $p\in\N$ be prime, $f\in\N$ and $q:=p^f\ge 3$ and suppose that $G=\PSL_3^\epsilon(q)$.
	Moreover let $d:=(3,q-\epsilon),\alpha:=\frac{q^2+\epsilon q+1}{d}$ and suppose that $l:=[G,g,1 \mid [\alpha,1]]$ is a Hurwitz datum. Then $l \in \Ri_4^*$.
\end{lemma}

\prf{
As $\alpha$ and $\beta:=\frac{(q^2-1)}{d}$ are coprime to $p$, every element of $G$ of order $\alpha$ or $\beta$ is semisimple. 
We also note that $\alpha$ and $\frac{|G|}{\alpha}$ are coprime.

Let $K$ be a conjugacy class of $G$ of elements of order $\beta$ and with centraliser of order $\beta$ in $G$. (For more details see Table 2 in \cite{SimFra73}, penultimate line.)
Next we let $c \in G$ be of order $\alpha$. With the corollary after Theorem 2 in \cite{Gow} we find
$a \in K$ and $b \in G$ such that $c^{-1}=[a,b]$. 
We claim that $\{a,b,c\}$ is a generating set for $G$ that satisfies all conditions from Lemma \ref{crit}.
In fact, the only property that is left to prove is (c).
So we let $U:=\langle a,b,c\rangle$.
As $a,c\in U$, it follows that $\abs U$ is divisible by the lowest common multiple of $\alpha$ and $\beta$, hence by $\alpha \cdot \beta$. Assume for a contradiction that $G$ has a maximal subgroup
$M$ that contains $U$. Then we
inspect  Theorem 6.5.3 in \cite{GLS3}, bearing in mind that $\alpha \cdot \beta$ divides $|M|$.

Case (a) is impossible because $\alpha$ divides $|M|$. (For more details see for example
3.3.3 and 3.6.2 in \cite[p. 47, p. 67]{Wilson} as well as \cite[p. 19, p. 118]{Taylor}.)

The possibilities
$\PSL_3(q_0)$, $\PSU_3(q_0)$, $\PGL_3(q_0)$ or $\PGU_3(q_0)$ (where $q_0$ is a proper divisor of $q$)
do not occur because of their orders.

Moreover $|M| \notin \{\abs{\PSL_2(q)}, d^{-1}3(q-\epsilon)^2,3\alpha\}$, because $q\ge 3$, and this excludes the possibilities (b), (c) and (e), (f), (g) of Theorem 6.5.3 in \cite{GLS3}.
As $\alpha$ is coprime to 6, we see that $M \not \cong\PSU_3(2)$ and $M \not \cong\PGU_3(2)$.\\
We look at the cases (d), (h), (i), (j), and (k) from the theorem individually:\\
If $M\cong\PSL_3(2)$, then $\alpha=7$ and $q(q+\epsilon)=6$. As $q\ge 3$, this forces $q=3$ and $\epsilon=-1$.
Consequently $o(y)=\beta=8$, but $y\in M$ and  $\PSL_3(2)$ does not have an element of order $8$.\\
If $M \cong\Alt_6$ or $M \cong \M_{10}$, then the fact that $\alpha$ is coprime to $6$ implies that $\alpha=5$ and $q\cdot(q+\epsilon)=4$. This contradicts the fact that $q\ge 3$.\\
If $M \cong\Alt_7$, then $\alpha\in\{5,7,5\cdot 7\}$. We have already excluded the case $\alpha \neq 5$, just above. In the paragraph before we saw that $\alpha=7$ means that $q=3$, $\epsilon=-1$ and $\beta=8$. However, $\Alt_7$ does not have any elements of order 8. Then $\alpha=35$, but $\Alt_7$ does not have any elements of order $35$ either.\\

Therefore $U=G$ and Lemma \ref{crit} yields that $l$ has a witness.
The branching information and Theorem \ref{list3} imply that all non-trivial elements of $G$ have three fixed points or none, independently of the choice of the witness, so $l \in \Ri_4^*$.
}

\textbf{Proof of Theorem 1.1:}

\prf{Suppose that $G$ is as described in the first hypothesis of the theorem and let $\Omega$ be a $G$-orbit of $X$ on which $G$ acts with fixity 3. Then Theorem \ref{list3} gives the possibilities for $G$, the size of $\Omega$ and the orders of the point stabilisers. In each case there is only one possibility for a fixity 3 action, and every non-trivial element in a point stabiliser
fixes exactly three points. Then the hypothesis that $G$ acts with fixity at most 4 in total implies that there is only one orbit on which $G$ acts as it does on $\Omega$.
If $\Omega$ is the only non-regular orbit, then $G$ acts with branching datum as in Lines 1--8 of Table\ref{DataList}.\\
However, it is possible for $G$ to have further non-regular orbits on which it acts with fixity at most 2.
Since $G$ is not a Frobenius group, we only need to consider fixity 2 actions.
So we inspect Theorem \ref{list3} and compare it with Theorem 1.2 in \cite{MW} for groups that allow for mixed fixity action, and then we check Lemma 3.11 and Lemma 3.13 of \cite{MW} for the specific possibilities and the orders of the point stabilisers.
If we also take into account that $G$ acts with global fixity at most 4, then the number of non-regular orbits with each given action is restricted, and this leads to the possibilities mentioned in Remark \ref{mixed} and hence to the specific branching data in the table:\\
If $G=\PSL_2(7)$, then there can be one or two orbits where $G$ acts with point stabilisers of order 3, and
also one or two orbits where $G$ acts with point stabilisers of order 4. Together with
the orbit on which $G$ acts with fixity 3 this gives the possible branching data in Lines 9--16
of Table \ref{DataList}.\\
If $G=\PSL_3(4)$, then the only fixity 2 action occurs with point stabilisers of order $5$, which gives one or two possible orbits on which $G$ acts with fixity 2.
This gives the remaining lists in Table \ref{DataList}.\\

For the second statement we suppose that $l$ is one of the lists in the table and that it is a Hurwitz datum.
Then we replace $l$ by a list $l_0$ that contains the smallest possible cogenus and otherwise has the same entries as $l$.
For Line 7 of the table, this means that Lemma \ref{RHD:LU3} applies, so there is nothing left to prove.
For all the other lines, we have seen how to check the list $l_0$, i.e. the list with the minimal possible cogenus, for the existence of a witness by using \texttt{GAP}.
Then Lemma 3.6 in \cite{SW} yields that there is a witness for $l$, as well, whenever $g,g_0$ are such that the list is a Hurwitz datum in the first place.
As the arguments for the action of $G$ on a witness $X$, and hence for the containment of $l$ in $\Ri_4^*$, are independent of $g$ and $g_0$, this concludes the proof.
}

In some sense, this theorem is the strongest result that could have been expected: all potential branching
data actually occur, and for all witnesses the group acts as described with global fixity at most 4.
This is different from the situation in \cite{SW} where for $\Alt_5$ we found a Hurwitz datum without a witness, and it is also different from the situation for fixity 4, where for the groups $\PSL_2(7)$ and $\PSL_2(8)$ we found Hurwitz data without witnesses.
But this will be discussed in the next article!

	
\end{document}